\documentclass [11pt,twoside,b5paper]{article}
\usepackage{amsmath}
\usepackage{latexsym}
\usepackage{amssymb}
\usepackage{graphicx}
\usepackage{amsmath}
\usepackage{amssymb}
\usepackage{graphicx}
\usepackage{amsmath}
\usepackage{amsfonts}
\usepackage{euscript}
\usepackage{epsfig}
 \usepackage{color}
\usepackage{amsmath}
\usepackage{latexsym}
\usepackage{amssymb}
\usepackage{cite}
\usepackage{indentfirst}
\usepackage{enumerate}
\usepackage{tikz}

\def\enddemo{\qed \endtrivlist}
\expandafter\let\csname enddemo*\endcsname=\enddemo

\def\qedsymbol{\ifmmode\bgroup\else$\bgroup\aftergroup$\fi
  \vcenter{\hrule\hbox{\vrule height.6em\kern.6em\vrule}\hrule}\egroup}
\def\qed{\ifmmode\else\unskip\nobreak\fi\quad\qedsymbol}

 \newtheorem{thm}{Theorem}[section]
 \newtheorem{prop}{Proposition}[section]
 \newtheorem{lem}{Lemma}[section]

 \newcounter{slicica}[section]

\title{\bf Projective Curvature Tensors of Second Type Almost Geodesic Mappings}
\author{Nenad O. Vesi\'c}
\date{}

\usepackage[T1]{fontenc}

\newcommand\blfootnote[1]{%
  \begingroup
  \renewcommand\thefootnote{}\footnote{#1}%
  \addtocounter{footnote}{-1}%
  \endgroup
}

\begin{document}
  \maketitle

  \begin{abstract}
    We consider equitorsion second type almost geodesic mappings of
    a non-symmetric affine connection space in this article. Using
    different computational methods, we obtained some invariants of  these
    mappings. Last generalized Thomas projective parameter and Weyl projective
    tensor as invariants of a second type almost geodesic mapping
    of a non-symmetric affine connection space are further generalized here.

    \blfootnote{2010 \emph{Math. Subj. Classification}: Primary:
    53A55; Secondary: 53B05, 53C15, 53C22}
    \blfootnote{This research is
supported by Ministry of Education, Science and Technological
Development, Republic of Serbia, Grant No. 174012}
  \end{abstract}

  \section{Introduction}

A lot of research papers and monographs are dedicated to
developments of the theory of differential geometry
\cite{eisenhart1, ja, jamica, mica1, mica2, mica3, micasg,
micazlatja, minc2, minc3, pi1sim, pi2esim, pi2sim, rec1, rec2, rec3,
sinjukov} and its applications \cite{e1, e2, e3, zlativanov}.
Einstein (see \cite{e1, e2, e3}) concluded the symmetric affine
connection theory covers researches about a gravitation. The theory
electromagnetism is covered by anti-symmetric parts of affine
connections. The research about non-symmetric affine connected
spaces is started by L. P. Eisenhart \cite{eisenhart1}.

An $N$-dimensional manifold ${\cal M}_N$ endowed with a
non-symmetric affine connection $\nabla$ (affine connection
coefficients $L^i_{jk}$ and $L^i_{kj}$ are different) is said to be
\emph{the non-symmetric affine connection space} $\mathbb{GA}_N$.
Because of the previous mentioned non-symmetry of affine connection
coefficients it exists symmetric and anti-symmetric part of these
coefficients respectively defined as:

\begin{eqnarray}
  \widetilde S^i_{jk}=\frac12(L^i_{jk}+L^i_{kj})&\mbox{and}&
  \widetilde T^i_{jk}=\frac12(L^i_{jk}-L^i_{kj}).
  \label{eq:ST}
\end{eqnarray}

A symmetrization and an anti-symmetrization without division by
indices $i$ and $j$ will be denoted as $(i\ldots j)$ and $[i\ldots
j]$ respectively.

 The magnitude $\widetilde T^i_{jk}$ is a torsion
tensor of the space $\mathbb{GA}_N$. An affine connection space
$\mathbb A_N$ endowed with an affine connection $S$ which
coefficients coincide with the symmetric part $\widetilde S^i_{jk}$
of the affine connection coefficients $L^i_{jk}$ of the space
$\mathbb{GA}_N$ is said to be \emph{the associated space} of the
space $\mathbb{GA}_N$.\vspace{.2cm}

There are a lot of researchers interested for a development of  the
non-symmetric affine connection space theory. Some significant
results in this subject are obtained into the papers \cite{minc2,
minc3, ja, jamica,  mica1, mica2, mica3, micasg, micazlatja}.

Four kinds of covariant differentiation (see \cite{minc2}) with
regard to an affine connection of a non-symmetric affine connection
space $\mathbb{GA}_N$ are defined as:

\begin{eqnarray}
  a^i_{j\underset1|k}=a^i_{j,k}+L^i_{\alpha
  k}a^\alpha_j-L^\alpha_{jk}a^i_\alpha,&&
  a^i_{j\underset2|k}=a^i_{j,k}+L^i_{k\alpha}a^\alpha_j-
  L^\alpha_{kj}a^i_\alpha,\label{eq:covdev1}\\
  a^i_{j\underset3|k}=a^i_{j,k}+L^i_{\alpha
  k}a^\alpha_j-L^\alpha_{kj}a^h_\alpha,&&
  a^i_{j\underset4|k}=a^i_{j,k}+L^i_{k\alpha}a^\alpha_j-
  L^\alpha_{jk}a^i_\alpha,\label{eq:covdev4}
\end{eqnarray}
for a partial derivative denoted by comma and an indexed magnitude
$a^i_j$.

All of these covariant derivatives become restricted to a covariant
derivative

\begin{equation}
  a^i_{j;k}=a^i_{j,k}+\widetilde S^i_{\alpha k}a^\alpha_j-
  \widetilde S^\alpha_{jk}a^i_\alpha,
  \label{eq:covdev}
\end{equation}
of the magnitude $a^i_j$ with regard to an affine connection of the
associated space $\mathbb A_N$ of the space $\mathbb{GA}_N$.

For this reason, 
%
it exists only one curvature tensor

\begin{equation} R^i_{jmn}=
\widetilde S^i_{jm;n}-\widetilde
S^i_{jn;m},\label{eq:r}\end{equation} of the associated space
$\mathbb{A}_N$.

\subsection{Almost geodesic mappings of a space $\mathbb{GA}_N$}

In an attempt to generalize the term of geodesics N. S. Sinyukov
(see \cite{sinjukov}) defined an almost geodesic line of a symmetric
 affine connection space $\mathbb{A}_N$.
Consequently, he defined a term of an almost geodesic mapping $f$
between symmetric affine connection spaces $\mathbb{A}_N$ and
$\mathbb{\overline A}_N$. Sinyukov noticed three types
$\pi_1,\pi_2,\pi_3$ of almost geodesic mappings between symmetric
affine connection spaces. His research has been directly developed
by many authors in a lot of papers \cite{rec1, pi1sim, rec2, rec3,
pi2sim, pi2esim}.

The Sinyukov's generalization of geodesics is primary developed for
the case of a generalized affine connection space $\mathbb{GA}_N$ in
\cite{mica1, mica2, mica3}. In this space it exists four kinds of
covariant differentiation but these covariant derivatives are
reduced onto first two ones (\ref{eq:covdev1}) for the case of any
contra-variant tensor.

For this reason, there are two kinds of almost geodesic lines of the
space $\mathbb{GA}_N$ \cite{mica1, mica2, mica3, micasg, micazlatja,
ja, jamica} defined as a curve $\ell=\ell(t)$ which tangential
vector\linebreak $\lambda^i=d\ell^i/dt\neq0$ satisfies the following
equations:

\begin{eqnarray}
  \underset\theta{\overline\lambda}{}^i_{(2)}=
  \underset\theta{\overline a}(t)\lambda^i+
  \underset\theta{\overline
  b}(t)\underset\theta{\overline\lambda}{}^i_{(1)},&
  \underset\theta{\overline\lambda}{}^i_{(1)}=\lambda^i_{\underset\theta\|\alpha}\lambda^\alpha,&
  \underset\theta{\overline\lambda}{}^i_{(2)}=
  \underset\theta{\overline\lambda}{}^i_{(1)\underset\theta\|\alpha}\lambda^\alpha,
  \label{eq:aglGAN}
\end{eqnarray}
$\theta=1,2$, for covariant differentiation of the $\theta$-th kind
with regard to affine connection of the space $\mathbb{G\overline
A}_N$ denoted by $\underset\theta{\|}$.\vspace{.2cm}

Because of two kinds of almost geodesic lines of this space
 a mapping
$f:\mathbb{GA}_N\rightarrow\mathbb{G\overline A}_N$ is \textbf{the
almost geodesic mapping of a $\theta$-th kind}, $\theta=1,2$, if any
geodesic line of the space $\mathbb{GA}_N$ it turns into an almost
geodesic line of the $\theta$-th kind of the space
$\mathbb{G\overline A}_N$. For this reason, there are three types of
almost geodesic mappings of the space $\mathbb{GA}_N$ and any of
these three types have two kinds. A class of almost geodesic
mappings of a $\tau$-th type, $\tau=1,2,3$, and of a $\theta$-th
kind $\theta=1,2$ of the space $\mathbb{GA}_N$ is denoted as
$\underset\theta\pi{}_\tau$.

Basic equations of a second type almost geodesic mapping\linebreak
$f:\mathbb{GA}_N\rightarrow\mathbb{G\overline A}_N$ of a $\theta$-th
kind, $\theta=1,2$, are \cite{mica2}:

\begin{align}
&\aligned
    \overline L^i_{jk}=L^i_{jk}+\psi_j\delta^i_k+\psi_k\delta^i_j+\sigma_jF^i_k+\sigma_kF^i_j+\xi^i_{jk},\\
\endaligned\label{eq:basicns1}\\
&\aligned    F^i_{j\underset\theta|k}&+F^i_{k\underset\theta|j}+
F^i_\alpha
    F^\alpha_j\sigma_k+F^i_\alpha
    F^\alpha_k\sigma_j+(-1)^\theta\big(\xi^i_{j\alpha}F^\alpha_k+\xi^i_{k\alpha}F^\alpha_j\big)\\&=
    \mu_jF^i_k+\mu_kF^i_j+\nu_j\delta^i_k+\nu_k\delta^i_j,
  \endaligned\label{eq:basicns2}
\end{align}
for covariant vectors $\mu_j,\nu_j$, an affinor $F^i_j$ and an
anti-symmetric tensor $\xi^i_{jk}$.

A second type almost geodesic mapping
$f:\mathbb{GA}_N\rightarrow\mathbb{G\overline A}_N$ of a $\theta$-th
kind, $\theta=1,2$, satisfies the property of reciprocity \big(it is
an element of the class $\underset\theta\pi{}_2(e)$\big) if it saves
the affinor $F^i_j$ and its inverse mapping is a second type almost
geodesic mapping of the $\theta$-th kind. An almost geodesic mapping
$f$ of the space $\mathbb{GA}_N$ satisfies the property of
reciprocity (see \cite{mica2}) if and only if the affinor $F^i_j$
satisfies a relation

\begin{equation}
  F^i_\alpha F^\alpha_j=e\delta^i_j,\quad\quad\quad e=0,\pm1.
  \label{eq:ffee}
\end{equation}

\section{Invariants of second type almost geodesic mappings}

The aim of this paper is to find some new invariants of almost
geodesic mappings of a second type which satisfy the property of
reciprocity. The results in this subject obtained until now are
about the theories of special subclasses of the classes
$\underset\theta\pi{}_2,\theta=1,2$.

Motivated by Sinyukov's results, it is obtained (see \cite{mica2})
magnitudes

\begin{align}
    &\aligned
  \underset1T{}^i_{jk}=\widetilde S{}^i_{{jk}}-\frac1{e-F^2}\Big(\big(
  F\widetilde S{}^\alpha_{{k\alpha}}-F^\alpha_k\widetilde S{}^\beta_{{\alpha\beta}}\big)F^i_j+
  \big(F\widetilde S{}^\alpha_{{j\alpha}}-F^\alpha_j\widetilde S{}^\beta_{{\alpha\beta}}\big)F^i_k\Big),
  \endaligned\\\displaybreak[0]
  &\aligned
  \underset2{\hat
  T}{}^i_{jk}&=T^i_{jk}+eF^i_\alpha\big(F^\alpha_{(j\underset1|k)}-\widetilde T{}^\alpha_{\beta(k}F^\beta_{j)}\big)\\&-
  \frac
  e{1+N}F^\beta_\alpha\Big(\big(F^\alpha_{\beta\underset1|j}-\widetilde T{}^\alpha_{\gamma(\beta}F^\gamma_{j)}\big)\delta^i_k+
  \big(F^\alpha_{\beta\underset1|k}-\widetilde T{}^\alpha_{\gamma(\beta}F^\gamma_{k)}\big)\delta^i_j\Big),
  \endaligned
\end{align}
for $F=F^\alpha_\alpha$, $e-F^2\neq0$ and Thomas projective
parameter $T^i_{jk}$ of the associated space $\mathbb A_N$ in the
expression of the invariant $\underset1T{}^i_{jk}$, are invariants
of a canonical second type almost geodesic mapping of the first kind
of the space $\mathbb{GA}_N$.

Moreover, Weyl projective tensor of the space $\mathbb{G\hat A}_N$
which affine connection coefficients are $\hat L^i_{jk}=L^i_{jk}+
eF^i_\alpha
F^\alpha_{(j\underset1|k)}-eT^\alpha_{\beta(j}F^\beta_{k)}F^i_\alpha$
is an invariant of the canonical second type almost geodesic mapping
$f$ of the first kind. The aim of our following research is to find
some other more general invariants of special second type almost
geodesic mappings of the space $\mathbb{GA}_N$.\vspace{.2cm}

Let a mapping $f:\mathbb{GA}_N\rightarrow\mathbb{G\overline A}_N$ be
an equitorsion second type almost geodesic mapping of a $\theta$-th
kind, $\theta=1,2$, which satisfies the property of reciprocity. The
composition (\ref{eq:ffee}) involved into the basic equation
(\ref{eq:basicns2}) together with using of the fact the mapping $f$
is an equitorsion one \linebreak($\xi^i_{jk}=0$) involved into the
both of basic equations (\ref{eq:basicns1}, \ref{eq:basicns2})
proves it is satisfied relations

\begin{align}
&\overline
L^i_{jk}=L^i_{jk}+\psi_j\delta^i_k+\psi_k\delta^i_j+\sigma_jF^i_k+\sigma_kF^i_j,\label{eq:basicns1e}
\\&
F^i_{j\underset\theta|k}+F^i_{k\underset\theta|j}=\mu_jF^i_k+\mu_kF^i_j+
  \big(\nu_j-e\sigma_j\big)\delta^i_k+\big(\nu_k-e\sigma_k\big)\delta^i_j.
  \label{eq:basicns2e}
\end{align}

It is proved a following proposition is satisfied in this way.

\begin{prop}
  Let $f:\mathbb{GA}_N\rightarrow\mathbb{G\overline A}_N$ be an
  equitorsion
  second type almost geodesic mapping of a $\theta$-th kind, $\theta=1,2$,
   which satisfies the property
  of reciprocity. The equations
  \emph{(\ref{eq:basicns1e}, \ref{eq:basicns2e}, \ref{eq:ffee})} are
  basic equations of this mapping.\qed
\end{prop}

Based on the fact the second type almost geodesic mapping $f$
satisfies the property of reciprocity the corresponding magnitudes
$\overline\psi_i,\overline\sigma_i,\overline F^i_j$ which determine
an inverse mapping $f^{-1}$ of the mapping $f$ are \cite{mica2}

\begin{eqnarray*}
  \overline\psi_i=-\psi_i,&\overline\sigma_i=-\sigma_i,&\overline
  F^i_j=F^i_j.
\end{eqnarray*}

After contracting the basic equation (\ref{eq:basicns1e}) by indices
$i$ and $k$ and using the fact it is satisfied a relation
$\sigma_j=\frac12(\sigma_j-\overline\sigma_j)$ we obtain it is
satisfied an equation

\begin{equation}
\aligned \psi_j\!=\!\frac1{N+1}\big(\overline
L^\alpha_{j\alpha}\!-\!
L^\alpha_{j\alpha}\big)\!+\!\frac1{2(N+1)}\Big[
\big(\overline\sigma_j\overline
F\!+\!\overline\sigma_\alpha\overline F^\alpha_j\big)\!-\!
\big(\sigma_jF\!+\!\sigma_\alpha F^\alpha_j\big)\Big],
\endaligned
\end{equation}
for $F=F^\alpha_\alpha$ as above.

Using the previous expression of the magnitude $\psi_j$ we conclude
the basic equation (\ref{eq:basicns1e}) has a form

\begin{equation}
\overline L^i_{jk}=L^i_{jk}+\overline\omega^i_{jk}-\omega^i_{jk},
\label{eq:basicns1ee}
\end{equation}
for

\begin{equation}
  \aligned
  \omega^i_{jk}=&-\frac12(\sigma_jF^i_k+\sigma_kF^i_j)+
  \frac1{N+1}\big(L^\alpha_{j\alpha}\delta^i_k+L^\alpha_{k\alpha}\delta^i_j\big)\\&
  +\frac1{2(N+1)}\Big(\big(\sigma_jF+\sigma_\alpha
  F^\alpha_j\big)\delta^i_k+
  \big(\sigma_kF+\sigma_\alpha
  F^\alpha_k\big)\delta^i_j\Big),
  \endaligned
  \label{eq:omega2}
\end{equation}
and the magnitude $\overline\psi_j$ defined in the same manner as a
function of the corresponding elements of the space
$\mathbb{G\overline A}_N$.

The equation (\ref{eq:basicns1ee}) proves it is satisfied an
equality

\begin{equation*}
  \underset2{\overline{\cal T}}{}^i_{jk}=\underset2{\cal
  T}{}^i_{jk},
\end{equation*}
for
\begin{eqnarray}
  \underset2{\cal T}{}^i_{jk}=L^i_{jk}-\omega_{jk}^i&\mbox{and}&
  \underset2{\overline{\cal T}}{}^i_{jk}=\overline
  L^i_{jk}-\overline\omega^i_{jk}.
  \label{eq:Thompi2}
\end{eqnarray}

It is proved a following lemma in this way.

\begin{lem}
  Let $f:\mathbb{GA}_N\rightarrow\mathbb{G\overline A}_N$ be an
  equitorsion almost geodesic mapping of a second type which
  satisfies the property of reciprocity. A magnitude
  $\underset2{\cal T}{}^i_{jk}$ defined in the first of the
  expressions \emph{(\ref{eq:Thompi2})} is an invariant of the
  mapping $f$.\qed
\end{lem}

The invariant $\underset2{\cal T}{}^i_{jk}$ of a second type almost
geodesic mapping\linebreak
$f:\mathbb{GA}_N\rightarrow\mathbb{G\overline A}_N$ which satisfies
the property of reciprocity is said to be \textbf{the
$\pi_2$-generalized Thomas projective parameter}.\vspace{.2cm}

Let us generalize Weyl projective tensor of an equitorsion almost
geodesic mapping $f:\mathbb{GA}_N\rightarrow\mathbb{G\overline A}_N$
of the first type which satisfies the property of reciprocity. This
generalization will be realized just for first type almost geodesic
mappings. A result for the case of second type almost geodesic
mappings may be obtained in the same manner.

First of all, we can observe the magnitude $\omega{}^i_{jk}$ defined
into the equation (\ref{eq:omega2}) is symmetric by indices $j$ and
$k$. For this reason, symmetric parts $\widetilde S{}^i_{jk}$ and
$\widetilde{\overline S}{}^i_{jk}$ of affine connection coefficients
$L^i_{jk}$ and $\overline L^i_{jk}$ of the spaces $\mathbb{GA}_N$
and $\mathbb{G\overline A}_N$ satisfy a relation

\begin{equation}
  \widetilde{\overline S}{}^i_{jk}=\widetilde
  S{}^i_{jk}+\overline\omega{}^i_{jk}-\omega{}^i_{jk},
  \label{eq:basicns1s}
\end{equation}
for the above defined magnitudes $\omega^i_{jk}$ and
$\overline\omega{}^i_{jk}$.\vspace{.2cm}

Using the covariant derivative of the first kind (\ref{eq:covdev1})
we obtain a curvature tensor $R^i_{jmn}$  of the associated space
$\mathbb A_N$ has a form:

\begin{align}
  &\aligned
  R^i_{jmn}&=\widetilde S{}^i_{jm\underset1|n}-
  \widetilde S{}^i_{jn\underset1|m}-\widetilde T{}^i_{\alpha n}
  \widetilde S^\alpha_{jm}-\widetilde T{}^\alpha_{jm}\widetilde
  S{}^i_{\alpha n}\\&+
  \widetilde T{}^\alpha_{jn}\widetilde S{}^i_{\alpha m} + \widetilde
  T{}^i_{\alpha m}\widetilde S{}^\alpha_{jn}+2
  \widetilde T^\alpha_{mn}\widetilde S{}^i_{j\alpha}.
  \endaligned\label{eq:r|1}
\end{align}

Motivated by this result we are going to find a rule of change of
the curvature tensor $R^i_{jmn}$ bellow. Let us involve following
substitutions:

\begin{eqnarray}
\underset2{\cal U}{}^i_{jk}=\frac12\big(\underset2{\cal T}{}^i_{jk}+
\underset2{\cal T}{}^i_{kj}\big)&\mbox{and}&
\underset2{\overline{\cal
U}}{}^i_{jk}=\frac12\big(\underset2{\overline{\cal T}}{}^i_{jk}+
\underset2{\overline{\cal T}}{}^i_{kj}\big), \label{eq:UUdfn}
\end{eqnarray}
for the above obtained invariant $\underset2{\cal
T}{}^i_{jk}=\underset2{\overline{\cal T}}{}^i_{jk}$ of the mapping
$f$.

From the equations (\ref{eq:Thompi2}) and (\ref{eq:UUdfn}) together
with the above mentioned symmetry of the magnitude $\omega^i_{jk}$
from the equation (\ref{eq:omega2}) by indices $j$ and $k$ we
conclude it is satisfied equalities

\begin{eqnarray}
  \underset2{\cal U}{}^i_{jk}=\widetilde S{}^i_{jk}-\omega{}^i_{jk}&\mbox{and}&
  \underset2{\overline{\cal U}}{}^i_{jk}=\widetilde{\overline S}{}^i_{jk}-\overline{\omega}{}^i_{jk}.
  \label{eq:USomega2}
\end{eqnarray}

It is easy to be obtained covariant derivatives

\begin{eqnarray}
\underset2{\cal
U}{}^i_{{jm}\underset1|n}=\frac12\big(\underset2{\cal
T}{}^i_{{jm}\underset1|n}+\underset2{\cal
T}{}^i_{{mj}\underset1|n}\big)& \mbox{and}&
\underset2{\overline{\cal U}}{}^i_{{jm}\underset1\|n}=\frac12\big(
\underset2{\overline{\cal T}}{}^i_{jm\underset1\|n}+
\underset2{\overline{\cal
T}}{}^i_{mj\underset1\|n}\big),\quad\quad\quad \label{eq:UU}
\end{eqnarray}
of the previous defined magnitudes $\underset2{\cal U}{}^i_{{jm}}$
and $\underset2{\overline{\cal U}}{}^i_{{jm}}$  satisfy a relation

\begin{equation}
  \aligned
  \underset2{\overline{\cal U}}{}^i_{{jm}\underset1\|n}&=
  \underset2{{\cal U}}{}^i_{{jm}\underset1|n}+\widetilde{\overline
  S}{}^i_{\alpha n}\underset2{\overline{\cal U}}{}^\alpha_{{jm}}-
  \widetilde{\overline S}{}^\alpha_{jn}\underset2{\overline{\cal
  U}}{}^i_{{\alpha m}}-\widetilde{\overline S}{}^\alpha_{mn}
  \underset2{\overline{\cal U}}{}^i_{{j\alpha}}
  \\&-\widetilde{S}{}^i_{\alpha n}\underset2{{\cal U}}{}^\alpha_{{jm}}+
  \widetilde{S}{}^\alpha_{jn}\underset2{{\cal
  U}}{}^i_{{\alpha m}}+\widetilde{S}{}^\alpha_{mn}
  \underset2{{\cal U}}{}^i_{{j\alpha}}.
  \endaligned\label{eq:S|S|}
\end{equation}

The equations (\ref{eq:USomega2}, \ref{eq:S|S|}) prove it is
satisfied a following proposition.

\begin{prop}
  Let $f:\mathbb{GA}_N\rightarrow\mathbb{G\overline A}_N$ be an
  equitorsion
  second type almost geodesic mapping of the first kind between
   non-symmetric affine connection
   spaces $\mathbb{GA}_N$ and $\mathbb{G\overline A}_N$.
   Covariant derivatives $\widetilde{S}{}^i_{jm\underset1|n}$ and
   $\widetilde{\overline S}{}^i_{jm\underset1\|n}$ of
   symmetric parts $\widetilde S{}^i_{jm}$ and $\widetilde{\overline S}{}^i_{jm}$
    of the corresponding affine connection coefficients $L^i_{jm}$
    and $\overline L^i_{jm}$ satisfy a relation

    \begin{equation}
      \aligned
      \widetilde{\overline
      S}{}^i_{jm\underset1\|n}&=\widetilde{S}{}^i_{jm\underset1|n}+
      \overline\omega{}^i_{jm\underset1\|n}-\omega{}^i_{jm\underset1|n}+
      \widetilde{\overline S}{}^i_{\alpha n}\underset2{\overline{\cal
      U}}{}^\alpha_{jm}\!-\!\widetilde{\overline
      S}{}^\alpha_{jn}\underset2{\overline{\cal U}}{}^i_{\alpha m}\!-\!
      \widetilde{\overline S}{}^\alpha_{mn}\underset2{\overline{\cal
      U}}{}^i_{j\alpha}\\&
      \!-\!\widetilde S{}^i_{\alpha n}\underset2{\cal U}{}^\alpha_{jm}\!+\!
      \widetilde S^\alpha_{jn}\underset2{\cal U}{}^i_{\alpha m}\!+\! \widetilde
      S{}^\alpha_{mn}\underset2{\cal U}{}^i_{j\alpha},
      \endaligned\label{eq:SS|}
    \end{equation}
    for the magnitudes $\underset2{\cal U}{}^i_{jk}$ and $\underset2{\overline{\cal U}}{}^i_{jk}$
     defined above.\qed
\end{prop}

Using the invariance $\widetilde{\overline
T}{}^i_{jk}=\widetilde{T}{}^i_{jk}$ and the consequent invariances
\linebreak$\underset2{\overline{\cal
T}}{}^\alpha_{jm}\widetilde{\overline T}{}^i_{\alpha
n}=\underset2{{\cal T}}{}^\alpha_{jm}\widetilde{T}{}^i_{\alpha n}$
such as $\underset2{\overline{\cal
T}}{}^i_{j\alpha}\widetilde{\overline T}{}^\alpha_{mn}=
\underset2{{\cal T}}{}^i_{j\alpha}\widetilde{T}{}^\alpha_{mn}$ we
obtain it is satisfied a following proposition.

\begin{prop}
  Let $f:\mathbb{GA}_N\rightarrow\mathbb{G\overline A}_N$ be an
  equitorsion
  second type almost geodesic mapping of the first kind between
   non-symmetric affine connection
   spaces $\mathbb{GA}_N$ and $\mathbb{G\overline A}_N$. Magnitudes
   $\widetilde S{}^\alpha_{jm}\widetilde T{}^i_{\alpha n},
   \widetilde S{}^i_{\alpha j}\widetilde T{}^\alpha_{mn}$ and its
   deformations
   $\widetilde{\overline S}{}^\alpha_{jm}\widetilde{\overline T}{}^i_{\alpha n},
   \widetilde{\overline S}{}^i_{\alpha j}\widetilde{\overline
   T}{}^\alpha_{mn}$ under the mapping $f$ satisfy equations

   \begin{align}
     &\widetilde{\overline S}{}^\alpha_{jm}\widetilde{\overline T}{}^i_{\alpha
     n}=\widetilde{S}{}^\alpha_{jm}\widetilde{T}{}^i_{\alpha
     n}+\overline\omega{}^\alpha_{jm}\widetilde{\overline
     T}{}^i_{\alpha n}-
     \omega{}^\alpha_{jm}\widetilde{T}{}^i_{\alpha
     n},\label{eq:Tcov}\\\displaybreak[0]&
     \widetilde{\overline S}{}^i_{j\alpha}\widetilde{\overline
     T}{}^\alpha_{mn}=
     \widetilde{S}{}^i_{j\alpha}\widetilde{T}{}^\alpha_{mn}+
     \overline\omega{}^i_{j\alpha}\widetilde{\overline
     T}{}^\alpha_{mn}-\omega{}^i_{j\alpha}\widetilde{T}{}^\alpha_{mn},\label{eq:Scov}
   \end{align}
   for the magnitude $\omega^i_{jk}$ defined into the equation
   \emph{(\ref{eq:omega2})}  and the corresponding one $\overline\omega{}^i_{jk}$.\qed
\end{prop}

If we contract the basic equation (\ref{eq:basicns2e}) by the
indices $i$ and $k$ we conclude it is satisfied a relation

\begin{equation}
  F_{\underset\theta|j}=\mu_jF+\mu_\alpha
  F^\alpha_j+(N+1)(\nu_j-e\sigma_j)-F^\alpha_{j\underset\theta|\alpha},
  \label{eq:Ftheta}
\end{equation}
$\theta=1,2$. This equation proves it is satisfied a following
proposition.

\begin{prop}
  Let $f:\mathbb{GA}_N\rightarrow\mathbb{G\overline A}_N$ be an
  equitorsion
  second type almost geodesic mapping of the first kind between
   non-symmetric affine connection
   spaces $\mathbb{GA}_N$ and $\mathbb{G\overline A}_N$ which satisfies the property of reciprocity.
   A covariant
   derivative $\omega^i_{jm\underset1|n}$ of the magnitude
   $\omega^i_{jk}$ defined into the equation
   \emph{(\ref{eq:omega2})} satisfies an equality

   \begin{equation}
     \aligned
     \omega^i_{jm\underset1|n}&=\frac1{N+1}\big(L^\alpha_{j\alpha\underset1|n}\delta^i_m+
     L^\alpha_{m\alpha\underset1|n}\delta^i_j\big)+\frac12(\nu_n-e\sigma_n)(\sigma_j\delta^i_m+\sigma_m\delta^i_j)\\&-\frac12\big(\sigma_{j\underset1|n}F^i_m+
     \sigma_{m\underset1|n}F^i_j+\sigma_jF^i_{m\underset1|n}+\sigma_mF^i_{j\underset1|n}\big)\\&+
     \frac1{2(N+1)}\big(\sigma_{j\underset1|n}F+\sigma_{\alpha\underset1|n}F^\alpha_j+\sigma_\alpha
     F^\alpha_{j\underset1|n}\big)\delta^i_m\\&+
     \frac1{2(N+1)}\big(\sigma_{m\underset1|n}F+\sigma_{\alpha\underset1|n}F^\alpha_m+\sigma_\alpha
     F^\alpha_{m\underset1|n}\big)\delta^i_j\\&+\frac1{2(N+1)}\big(
     \mu_nF+\mu_\alpha
     F^\alpha_n-F^\alpha_{n\underset1|\alpha}\big)(\sigma_j\delta^i_m+\sigma_m\delta^i_j),
     \endaligned\label{eq:omega|}
   \end{equation}
   for magnitudes $\mu_i,\nu_i$ used into the basic equation
   \emph{(\ref{eq:basicns2e})}.\qed
\end{prop}

  A difference
  $\widehat\Delta{}^i_{jmn}=\overline\omega{}^i_{jm\underset1\|n}-
  \omega^i_{jm\underset1|n}$, of the magnitudes
  $\overline\omega{}^i_{jm\underset1\|n}$ and $\omega^i_{jm\underset1|n}$
   satisfies a relation

  \begin{equation}
  \aligned
    \widehat\Delta{}^i_{jmn}&=\frac1{N+1}\Big(\big(\widetilde{\overline
    S}{}^\alpha_{j\alpha\underset1\|n}-\widetilde
    S^\alpha_{j\alpha\underset1|n}\big)\delta^i_m+\big(\widetilde{\overline
    S}{}^\alpha_{m\alpha\underset1\|n}-\widetilde
    S{}^\alpha_{m\alpha\underset1|n}\big)\delta^i_j\Big)\\&+
    \hat{\overline\rho}{}^i_{jmn}-\hat\rho{}^i_{jmn},
  \endaligned\label{eq:omega-omega}
  \end{equation}
  for

  \begin{align}
    &\aligned
    2\hat{\rho}{}^i_{jmn}&=-2L^\beta_{jn}\widetilde{T}{}^\alpha_{\beta\alpha}\delta^i_m
    -2L^\beta_{mn}\widetilde{T}{}^\alpha_{\beta\alpha}\delta^i_j+(\nu_n-e\sigma_n)(\sigma_j\delta^i_m+\sigma_m\delta^i_j)\\&-
    (\sigma_{j\underset1|n}F^i_m+
    \sigma_{m\underset1|n}F^i_j+
    \sigma_j
    F^i_{m\underset1|n}+\sigma_m
    F^i_{j\underset1|n})\\&+\frac1{N+1}\big(\sigma_{j\underset1|n}
    F+\sigma_{\alpha\underset1|n}F^\alpha_j+
    \sigma_\alpha
    F^\alpha_{j\underset1|n}\big)\delta^i_m\\&+
    \frac1{N+1}\big(\sigma_{m\underset1|n} F+
    \sigma_{\alpha\underset1|n}F^\alpha_m+
    \sigma_\alpha
    F^\alpha_{m\underset1|n}\big)\delta^i_j\\&+\frac1{N+1}
    \big(\mu_nF+\mu_\alpha
    F^\alpha_n-
    F^\alpha_{n\underset1|\alpha}\big)(\sigma_j\delta^i_m+\sigma_m\delta^i_j),
    \endaligned\label{eq:rho}
  \end{align}
  a magnitude $\hat{\overline\rho}^i_{jmn}$ from the space $\mathbb{G\overline A}_N$
  analogue to the magnitude $\hat\rho{}^i_{jmn}$,
   the magnitudes $\mu_i,\nu_i$ from the equation
  {(\ref{eq:basicns2e})} and the corresponding ones
  $\overline\mu_i,\overline\nu_i$.\vspace{.2cm}

  From the  equations
  (\ref{eq:omega-omega}, \ref{eq:rho}) we conclude it exists a
  magnitude

  \begin{equation}
    \aligned
    \underset1{\hat\upsilon}{}_{ij}&=-\frac1{N+1}\widetilde{\overline
    S}{}^\alpha_{i\alpha\underset1\|j}
    +\overline L^\beta_{ij}\widetilde{\overline
    T}{}^\alpha_{\beta\alpha}-(\overline\nu_j-e\overline\sigma_j)\sigma_i\\&+\frac1{2(N+1)}
    \big(\overline\sigma_{i\underset1\|j}\overline F+
    \overline\sigma_{\alpha\underset1\|j}\overline F^\alpha_i+
    \overline\sigma_\alpha\overline
    F^\alpha_{i\underset1\|j}\big)\\&-\frac1{2(N+1)}\big(\overline\mu_j\overline
    F+\overline\mu_\alpha\overline F^\alpha_j-\overline
    F^\alpha_{j\underset1\|\alpha}\big)\overline\sigma_i\\&+
    \frac1{N+1}\widetilde S{}^\alpha_{i\alpha\underset1|j}-
    L^\beta_{ij}\widetilde
    T^\alpha_{\beta\alpha}+(\nu_j-e\sigma_j)\sigma_i\\&-
    \frac1{2(N+1)}\big(\sigma_{i\underset1|j}F+\sigma_{\alpha\underset1|j}F^\alpha_i+
    \sigma_\alpha F^\alpha_{i\underset1|j}\big)\\&+
    \frac1{2(N+1)}\big(\mu_jF+\mu_\alpha F^\alpha_j-
    F^\alpha_{j\underset1|\alpha}\big)\sigma_i,
    \endaligned
  \end{equation}
such that the equation (\ref{eq:SS|}) has a form

\begin{eqnarray}
  \aligned
  \widetilde{\overline S}{}^i_{jm\underset1\|n}&=\widetilde
  S{}^i_{jm\underset1|n}-\delta^i_m\underset1{\hat\upsilon}{}_{jn}-\delta^i_j\underset1{\hat\upsilon}{}_{mn}\\&-\frac12
  \big(\overline\sigma_{j\underset1\|n}\overline
  F^i_m+\overline\sigma_{m\underset1\|n}\overline F^i_j+
  \overline\sigma_j\overline F^i_{m\underset1\|n}+\overline\sigma_m
  \overline F^i_{j\underset1\|n}\big)\\&+\frac12
  \big(\sigma_{j\underset1|n}
  F^i_m+\sigma_{m\underset1|n}F^i_j+
  \sigma_jF^i_{m\underset1|n}+\sigma_m
  F^i_{j\underset1|n}\big)\\&+
      \widetilde{\overline S}{}^i_{\alpha n}\underset2{\overline{\cal
      U}}{}^\alpha_{jm}\!-\!\widetilde{\overline
      S}{}^\alpha_{jn}\underset2{\overline{\cal U}}{}^i_{\alpha m}\!-\!
      \widetilde{\overline S}{}^\alpha_{mn}\underset2{\overline{\cal
      U}}{}^i_{j\alpha}
      \!-\!\widetilde S{}^i_{\alpha n}\underset2{\cal U}{}^\alpha_{jm}\!+\!
      \widetilde S^\alpha_{jn}\underset2{\cal U}{}^i_{\alpha m}\!+\! \widetilde
      S{}^\alpha_{mn}\underset2{\cal U}{}^i_{j\alpha}.
  \endaligned\label{eq:SS||}
\end{eqnarray}

Using the equations (\ref{eq:UUdfn}, \ref{eq:Tcov}, \ref{eq:Scov},
\ref{eq:SS||}) and the expression (\ref{eq:r|1}) of the curvature
tensors $R^i_{jmn}$ and $\overline R^i_{jmn}$ of the associated
spaces $\mathbb A_N$ and $\mathbb{\overline A}_N$ we obtain it is
satisfied the following equation

\begin{eqnarray}
  \aligned
  \overline
  R^i_{jmn}&=R^i_{jmn}-\delta^i_m\underset1{\hat\upsilon}{}_{jn}+\delta^i_n\underset1{\hat\upsilon}{}_{jm}-\delta^i_j\underset1{\hat\upsilon}{}_{[mn]}\\&
  -\frac12
  \big(\overline\sigma_{j\underset1\|n}\overline
  F^i_m+\overline\sigma_{m\underset1\|n}\overline F^i_j+
  \overline\sigma_j\overline F^i_{m\underset1\|n}+\overline\sigma_m
  \overline F^i_{j\underset1\|n}\big)\\&+\frac12
  \big(\sigma_{j\underset1|n}
  F^i_m+\sigma_{m\underset1|n}F^i_j+
  \sigma_jF^i_{m\underset1|n}+\sigma_m
  F^i_{j\underset1|n}\big)\\&
  +\frac12
  \big(\overline\sigma_{j\underset1\|m}\overline
  F^i_n+\overline\sigma_{n\underset1\|m}\overline F^i_j+
  \overline\sigma_j\overline F^i_{n\underset1\|m}+\overline\sigma_n
  \overline F^i_{j\underset1\|m}\big)\\&-\frac12
  \big(\sigma_{j\underset1|m}
  F^i_n+\sigma_{n\underset1|m}F^i_j+
  \sigma_jF^i_{n\underset1|m}+\sigma_n
  F^i_{j\underset1|m}\big)
  \\&+
      2\widetilde{\overline S}{}^i_{\alpha n}\widetilde{\overline{S}}{}^\alpha_{jm}\!-\!
      2\widetilde{\overline
      S}{}^\alpha_{jn}\widetilde{\overline{S}}{}^i_{\alpha m}
      \!-\!2\widetilde S{}^i_{\alpha n}\widetilde{S}{}^\alpha_{jm}\!+\!
      2\widetilde S^\alpha_{jn}\widetilde{S}{}^i_{\alpha m}
      \\&-
      \overline\omega{}^\alpha_{jm}{\overline
      L}{}^i_{\alpha n}+\omega^\alpha_{jm}L{}^i_{\alpha
      n}-\overline\omega{}^i_{\alpha n}{\overline
      L}{}^\alpha_{jm}+\omega^i_{\alpha n}L{}^\alpha_{jm}\\&+\overline\omega{}^i_{\alpha
      m}{\overline L}{}^\alpha_{jn}\!-\!\omega^i_{\alpha
      m}L{}^\alpha_{jn}\!+\!\overline\omega{}^\alpha_{jn}{\overline
      L}{}^i_{\alpha m}\!-\!\omega^\alpha_{jn}L{}^i_{\alpha
      m}\!+\!2\overline\omega{}^i_{j\alpha}\widetilde{\overline
      T}{}^\alpha_{mn}\!-\!2\omega^i_{j\alpha}\widetilde
      T^\alpha_{mn}\\&=
      R^i_{jmn}-\delta^i_m\underset1{\hat\upsilon}{}_{jn}+\delta^i_n\underset1{\hat\upsilon}{}_{jm}-
      \delta^i_j\underset1{\hat\upsilon}{}_{[mn]}+\underset1{\cal F}{}^i_{jmn}-
      \underset1{\overline{\cal F}}{}^i_{jmn},
  \endaligned\label{eq:rr}
\end{eqnarray}
for

\begin{equation}
  \aligned
  \underset1{\cal F}{}^i_{jmn}&=\frac12
  \big(\sigma_{j\underset1|n}
  F^i_m+\sigma_{m\underset1|n}F^i_j+
  \sigma_jF^i_{m\underset1|n}+\sigma_m
  F^i_{j\underset1|n}\big)\\&-\frac12
  \big(\sigma_{j\underset1|m}
  F^i_n+\sigma_{n\underset1|m}F^i_j+
  \sigma_jF^i_{n\underset1|m}+\sigma_n
  F^i_{j\underset1|m}\big)\\&\!-\!2\widetilde S{}^i_{\alpha n}\widetilde{S}{}^\alpha_{jm}\!+\!
      2\widetilde S^\alpha_{jn}\widetilde{S}{}^i_{\alpha
      m}\\&+\omega^\alpha_{jm}L^i_{\alpha n}\!+\!\omega^i_{\alpha
      n}L^\alpha_{jm}\!-\!\omega^i_{\alpha m}L^\alpha_{jn}\!-\!
      \omega^\alpha_{jn}L^i_{\alpha
      m}\!-\!2\omega^i_{j\alpha}\widetilde T^\alpha_{mn},
  \endaligned\label{eq:f1}
\end{equation}
the magnitude $\omega^i_{jk}$ defined into the equation
(\ref{eq:omega2}) and the corresponding magnitude
$\underset1{\overline{\cal F}}{}^i_{jmn}$.

After contracting the equation (\ref{eq:rr}) by indices $i$ and $n$
we conclude Ricci tensors $R_{jm}$ and $\overline R_{jm}$ of the
associated spaces $\mathbb A_N$ and $\mathbb{\overline A}_N$ satisfy
a relation

\begin{equation}
  \overline
  R_{jm}=R_{jm}+(N-1)\underset1{\hat\upsilon}{}_{jm}+\underset1{\hat\upsilon}{}_{[jm]}+
  \underset1{\cal F}{}_{jm}-\underset1{\overline{\cal F}}{}_{jm},
  \label{eq:ricric}
\end{equation}
for $\underset1{\cal F}{}_{jm}=\underset1{\cal
F}{}^\alpha_{jm\alpha}$ and $\underset1{\overline{\cal
F}}{}_{jm}=\underset1{\overline{\cal F}}{}^\alpha_{jm\alpha}$.

After alternating the equation (\ref{eq:ricric}) by indices $j$ and
$m$ we conclude it is satisfied a relation

\begin{equation}
  (N+1)\underset1{\hat\upsilon}{}_{[jm]}=\overline R_{[jm]}-R_{[jm]}-
  \underset1{\cal F}{}_{[jm]}+\underset1{\overline{\cal
  F}}{}_{[jm]}.
  \label{eq:upsilon[]}
\end{equation}

From the equations (\ref{eq:ricric}) and (\ref{eq:upsilon[]}) we
conclude it is satisfied an expression

\begin{equation}
  \aligned
  (N^2-1)\underset1{\hat\upsilon}{}_{jm}&=(N\overline
  R_{jm}+\overline R_{mj})-(NR_{jm}+R_{mj})\\&+
  (N\underset1{\overline{\cal F}}{}_{jm}+
  \underset1{\overline{\cal F}}{}_{mj})-
  (N\underset1{{\cal F}}{}_{jm}+
  \underset1{\cal F}{}_{mj}).
  \endaligned\label{eq:upsilon}
\end{equation}

After involving the results (\ref{eq:upsilon[]}, \ref{eq:upsilon})
in the equation (\ref{eq:rr}) we obtain it is satisfied an equality

\begin{equation*}
  \underset2{{\overline{\cal W}}}{}^i_{jmn}=\underset2{\cal
  W}{}^i_{jmn},
\end{equation*}
where we denoted by

\begin{equation}
  \aligned
  \underset2{\cal
  W}{}^i_{jmn}&=R^i_{jmn}\!+\!\frac1{N+1}\delta^i_jR_{[mn]}\!+\!\frac
  N{N^2-1}\delta^i_{[m}R_{jn]}\!+\!\frac1{N^2-1}\delta^i_{[m}R_{n]j}\\&+
  \underset1{\cal
  F}{}^i_{jmn}\!+\!\frac1{N+1}\delta^i_j\underset1{\cal F}{}_{[mn]}+
  \frac N{N^2-1}\delta^i_{[m}\underset1{\cal F}{}_{jn]}+
  \frac1{N^2-1}\delta^i_{[m}\underset1{\cal F}{}_{n]j},
  \endaligned\label{eq:PI21inv}
\end{equation}
a geometric object of the space $\mathbb{GA}_N$, where the magnitude
$\underset1{\cal F}{}^i_{jmn}$ is defined into the equation
(\ref{eq:f1}) and for the corresponding magnitude $\underset1{\cal
F}{}_{ij}=\underset1{\cal F}{}^\alpha_{ij\alpha}$. The corresponding
magnitude $\underset2{\overline{\cal W}}{}^i_{jmn}$ of the space
$\mathbb{G\overline A}_N$ is defined in the same manner.

It is proved a following theorem is satisfied in this way.

\begin{thm}
  Let $f:\mathbb{GA}_N\rightarrow\mathbb{G\overline A}_N$ be an
  equitorsion second type almost geodesic mapping of the
  first kind which satisfies the property of reciprocity. The
  magnitude $\underset2{\cal W}{}^i_{jmn}$ defined into the equation
  \emph{(\ref{eq:PI21inv})} is an invariant of this mapping.\qed
\end{thm}

\section{Acknowledgements}

This paper is financially supported by the Serbian Ministry of
Education, Science and Technological Developments, Grant. No.
174012.

\noindent\textbf{Author's addresses:}

\begin{enumerate}
  \item[-] \emph{Nenad O. Vesi\'c};\\[3pt]
  \emph{Address}: Faculty of Sciences and
  Mathematics,\\ Vi\v segradska 33, 18000 Ni\v s,  Serbia;
  \\[1pt] \emph{e-mail}: vesko1985@pmf.ni.ac.rs
\end{enumerate}
\end{document}